\documentclass[11pt,reqno]{amsart}
\usepackage{amsmath}
\usepackage{amssymb}

\usepackage[all,2cell,ps]{xy}
\newcommand{\edge}[1]{\ar@{-}[#1]}

\usepackage{graphicx}
\usepackage{psfrag}
 \newtheorem{theorem}{Theorem}[section]
\newtheorem{prop}[theorem]{Proposition}

\newtheorem{cor}[theorem]{Corollary}

\theoremstyle{definition}

\newtheorem{defn}[theorem]{Definition}

\newtheorem{example}[theorem]{Example}

\newtheorem{conjecture}[theorem]{Conjecture}

\def\m{\mathcal M}
\def\n{\mathcal N}

\def\N{\mathbb N}
\def\R{\mathbb R}
\def\Z{\mathbb Z}

\def\phi{\varphi}

\def\d{\mathfrak{D}}
\def\ccc{\mathfrak{C}}

\address{Oliver Jenkinson; School of Mathematical Sciences, Queen
  Mary University of London, Mile End Road, London, E1 4NS, UK.
  }

\begin{document}

\newcounter{algnum}
\newcounter{step}
\newtheorem{alg}{Algorithm}

\newenvironment{algorithm}{\begin{alg}\end{alg}}

\title{Ergodic optimization in dynamical systems}

\author{Oliver Jenkinson}

\begin{abstract}
Ergodic optimization is the study of problems relating to 
maximizing orbits, 
maximizing invariant measures and maximum ergodic averages. 
An orbit of a dynamical system is called $f$-maximizing if 
the time average of the real-valued function 
$f$ along the orbit is larger than along all other orbits, and
an invariant probability measure is called 
$f$-maximizing if it gives 
$f$ a larger space average than does any other invariant probability measure.
In this survey we consider the main strands of ergodic optimization, beginning with an influential model problem, and the interpretation of ergodic optimization as the zero temperature limit of thermodynamic formalism.
We describe typical properties of maximizing measures for various spaces of functions, the key tool of adding a coboundary so as to reveal properties of these measures, as well as
certain classes of functions where the maximizing measure is known to
be Sturmian.
\end{abstract}

\maketitle


\section{Introduction}

For a real-valued function defined on the state space of
a dynamical system, the topic of \emph{ergodic optimization}
revolves around understanding its largest possible ergodic average.
Taking the dynamical system to be a map $T:X\to X$, 
and denoting the function by
$f:X\to\R$,
attention is focused on the supremum of time averages
$\lim_{n\to\infty} \frac{1}{n}\sum_{i=0}^{n-1} f(T^ix)$
over those $x\in X$ for which the limit 
exists,
or alternatively (and in nice cases equivalently) on the supremum of space averages
$\int f\, d\mu$ over probability measures $\mu$ which are invariant under $T$.

In the most classical setting of a topological dynamical system, with $X$ a compact metric space and $T:X\to X$ continuous, and if 
$f$ is continuous, then the above suprema coincide.
Indeed the common value is a maximum, 
as
the weak$^*$ compactness of the set $\m_T$ of $T$-invariant
Borel probability measures 
guarantees
some $m\in\m_T$ satisfying
\begin{equation}\label{maxmeasdef}
\int f\, dm = \max_{\mu\in\m_T} \int f\, d\mu =: \beta(f) \,,
\end{equation}
and there exists $x\in X$ with $\lim_{n\to\infty} \frac{1}{n}\sum_{i=0}^{n-1} f(T^ix)=\beta(f)$, since $m$ may be
taken to be ergodic 
and $x$ an $m$-generic point.
Any such $m\in\m_T$ is called a \emph{maximizing measure} for $f$, and $\beta(f)$ is the \emph{maximum ergodic average}.

Ergodic optimization originated in the 1990s, with much early work focused on fixing a specific map $T$ and studying the dependence of the maximizing measure on a function $f$ which varied in some finite dimensional space $V$. Indeed a certain model problem 
(see \S \ref{modelsection})
consisting of $T$ the doubling map on the circle,
and $V$ the 2-dimensional vector space 
of degree-1 trigonometric polynomials,
turned out to be influential: 
various subsequent results were suggested
either by the behaviour of this model, or by the techniques used to understand it.
In this model,
any non-zero function in $V$ has a unique maximizing measure, 
this measure
is usually periodic (i.e.~supported on a single periodic orbit), though not always periodic.
The natural occurrence of non-periodic maximizing measures was itself somewhat surprising
(and had ramifications
 in related 
areas
\cite{bouschjenkinson,bouschmairesse}), while the apparent rarity of non-periodic maximizing measures 
anticipated
the 
programme 
(described here in \S \ref{typicalsection})
of establishing 
analogous results for $V$ an infinite dimensional function space
(e.g.~the space of Lipschitz functions) 
and investigating 
further
 generic properties of maximizing measures
(see \S \ref{othertypicalsection}).

The specific maximizing measures arising in the model problem of \S \ref{modelsection}, so-called \emph{Sturmian} measures,
turned out to be unexpectedly ubiquitous in a variety of ergodic optimization problems
(which we describe in \S \ref{ergodicdominancesection}),
encompassing similar low-dimensional function spaces, 
certain infinite dimensional cones of functions,
 and problems concerning the joint spectral radius of matrix pairs.
Various ideas
used to resolve the model problem 
have been the subject of subsequent research; most notably, the prospect of adding a coboundary to $f$ so as to reveal
properties of its maximizing measure
has been the cornerstone of much recent work
(described in  \S \ref{revelationsection} and \S \ref{revelationtheoremssection}), 
with
many authors 
equally inspired by parallels with Lagrangian dynamical systems.

Another significant strand of research in ergodic optimization, again already present in early works,
was 
its
interpretation (see \S \ref{zerotemperaturesection})
as a limiting \emph{zero temperature} version of the more classical thermodynamic formalism, with maximizing measures 
(referred to as \emph{ground states} by physicists)
arising as
zero temperature accumulation points of equilibrium measures;
work in this area has primarily focused on understanding convergence and non-convergence in the zero temperature limit.

\section{Fundamentals}

Let $\d$ denote the set of pairs $(X,T)$ where $X=(X,d)$ is a compact metric space and $T:X\to X$ is continuous.
For $(X,T)\in\d$, the set $\m_T$ of $T$-invariant Borel probability measures
is compact when equipped with the weak$^*$ topology.

Let $\ccc$ denote the set of triples $(X,T,f)$, where $(X,T)\in\d$ and $f:X\to\R$ is continuous.
For $X$ a compact metric space, let $C(X)$ denote the set of continuous real-valued functions on $X$,
equipped with the supremum norm $\|f\|_\infty=\max_{x\in X}|f(x)|$. 
Let $Lip$ denote the set of Lipschitz real-valued functions on $X$, 
with $Lip(f):=\sup_{x\neq y} |f(x)-f(y)|/d(x,y)$, and Banach norm
$\|f\|_{Lip}=\|f\|_\infty+ Lip(f)$.

\begin{defn}
For $(X,T,f)\in\ccc$,
the quantity $\beta(f) = \beta(T,f) =  \beta(X,T,f)$ defined by
$$\beta(f) = \max_{\mu\in\m_T} \int f\, d\mu$$
is the \emph{maximum ergodic average}.
Any $m\in\m_T$ satisfying $\int f\, dm = \beta(f)$ is an \emph{$f$-maximizing measure},
and $\m_{\max}(f)=\m_{\max}(T,f)=\m_{\max}(X,T,f)$ denotes
the collection of
such measures.
\end{defn}

While we adopt the convention that optimization means maximization, occasional
mention will be made of
the \emph{minimum ergodic average} $$\alpha(f) = \min_{\mu\in\m_T} \int f\, d\mu = -\beta(-f)\,,$$
and the set 
$\m_{\min}(f)=\{m\in\m_T: \int f\, dm= \alpha(f)\}$
of  \emph{minimizing measures} for $f$. 
The closed interval $[\alpha(f),\beta(f)] = \{ \int f\, d\mu:\mu\in\m_T\}$ is the \emph{set of ergodic averages}\footnote{This set, and its 
generalisation
for $f$ taking values in higher dimensional spaces, is often referred to as the
\emph{rotation set} (see e.g.~\cite{blokh, garibaldilopes1, gellermisiurewicz, jenkinsoncras, jenkinsontams, kucherenkowolf1, kucherenkowolf2, ziemian}), while in the context of multifractal analysis it is sometimes referred to as the \emph{spectrum} of (Birkhoff) ergodic averages.}.

The maximum ergodic average admits a number of alternative characterisations involving time averages
(see e.g.~\cite[Prop~2.2]{jeo}):

\begin{prop}\label{alternative}
For $(X,T,f)\in\ccc$,
the maximum ergodic average $\beta(f)$ satisfies
\begin{equation}\label{fourfoldequality}
\beta(f)
=
\sup_{x\in X_{T,f}}
\lim_{n\to\infty} \frac{1}{n} S_n
  f(x) 
  =
\sup_{x\in X} \limsup_{n\to\infty} \frac{1}{n}
  S_nf(x)
  =
   \limsup_{n\to\infty} \frac{1}{n} \sup_{x\in X}
S_nf(x)
  \,,
  \end{equation}
  where $S_nf = \sum_{i=0}^{n-1} f\circ T^i$, and 
  $X_{T,f}=\{x\in X:  \lim_{n\to\infty} \frac{1}{n} S_nf(x)\text{ exists}\}$. 
\end{prop}

The following is well known (see e.g.~\cite[Prop.~2.4]{jeo}):

\begin{prop}
\label{mmaxf}
If $(X,T,f)\in\ccc$ then:

\item[$\,$ (i)] There exists at least one $f$-maximizing measure.
\item[$\,$ (ii)]
$\m_{\max}(f)$ is 
compact.
 \item[$\,$ (iii)]
$\m_{\max}(f)$ is a simplex, and in particular convex.
\item[$\,$ (iv)]
The extreme points of $\m_{\max}(f)$ are precisely those
$f$-maximizing measures which are ergodic.  In particular, there is at
least one ergodic $f$-maximizing measure.
\end{prop}

In \S \ref{typicalsection} and \S 8 we shall consider \emph{typical} properties of maximizing measures
in various spaces $V$ of real-valued functions on $X$.
The following result 
(see e.g.~\cite[Thm.~3.2]{jeo}, and in other forms see \cite{bousch2, contreraslopesthieullen, conzeguivarch})
guarantees that, for all of the function spaces $V$ considered,
\emph{uniqueness} of the maximizing measure is typical in $V$
(though clearly there exist $f\in V$ such that $\m_{\max}(f)$ is not a singleton, provided $\m_T$ is not a singleton,
most obviously $f\equiv 0$).


\begin{theorem}
\label{genericuniqueness} {\bf (Typical uniqueness of maximizing measures)}

\noindent
If $(X,T)\in\d$, and $V$
is a topological vector space which is densely and continuously
embedded in $C(X)$, then $\{f\in V: \m_{\max}(f)\text{ is a singleton}\}$ is a 
residual subset of $V$.
\end{theorem}

If $T$ and $f$ are continuous, but $X$ is non-compact, a number of difficulties potentially arise.
Assuming $\m_T$ is non-empty, we may define $\beta(f)=\sup_{\mu\in\m_T}\int f\, d\mu$, though in general there need not exist any maximizing measures, and any one of the equalities in (\ref{fourfoldequality}) may fail to hold
(see e.g.~\cite{jenkinsonmauldinurbanskinoncompact}). The most commonly studied example of a non-compact $X$
is a countable alphabet subshift of finite type, where a number of sufficient conditions have been given for existence of maximizing measures (see e.g.~\cite{bissacotfreire, bissacotgaribaldi, iommi, jenkinsonmauldinurbanskisft, jenkinsonmauldinurbanskinoncompact}), while \cite{davieurbanskizdunik} includes applications to 
(non-compact) Julia sets $X\subset\mathbb{C}$ for 
maps $T$ in the exponential family.

Note that versions of ergodic optimization 
have also been investigated in discrete time settings slightly different from the one described here, notably
the case where $\m_T$ is a singleton (see \cite{bremontbuczolich}),
in the context of non-conventional ergodic averages (see \cite{avramidou}),
or when the optimization is over a restricted subset of $\m_T$ (see \cite{zhaoobservable}).
Generalisations of ergodic optimization include optimal tracking for dynamical systems (see \cite{mcgoffnobel}),
and ergodic dominance (see \S \ref{ergodicdominancesection}).

\section{A model problem}\label{modelsection}

The map $T(x)= 2x \pmod 1$
on the circle $X=\R/\Z$ is a standard example of a hyperbolic dynamical system, and the functions
$f(x)=\cos 2\pi x$ and $g(x)=\sin 2\pi x$ are arguably the most natural non-constant functions on $X$.
While the $f$-maximizing measure is easily seen to be the Dirac measure at the fixed point 0, the $g$-maximizing measure is rather less obvious (it turns out to be the periodic measure on the orbit coded by 0001).
This standard choice of $T$, and the naturalness of $f$ and $g$, prompted 
several early authors to investigate
those $T$-invariant measures which are maximizing for functions in the 2-dimensional vector space $V$ spanned by $f$ and $g$. A rather complete understanding of this model problem has been provided by Bousch \cite{bousch1}, following earlier partial progress \cite{conzeguivarch, huntott1, huntott2, jenkinson1, jenkinson2, jenkinsonetds}, and the results
in this case already point to some more universal features of ergodic optimization.

While the space of degree-one trigonometric polynomials $V$ is 2-dimensional, the fact that a measure is maximizing for $v\in V$ if and only if it is maximizing for $cv$, where $c>0$, renders the problem a 1-dimensional one; specifically, to identify the maximizing measures for functions in $V$ it suffices to determine the maximizing measures for functions on the unit sphere in $V$, i.e.~those of the form
$v_\theta(x) = (\cos 2\pi\theta)f(x)+(\sin 2\pi\theta)g(x)=\cos 2\pi (x-\theta)$, for $\theta\in \R/\Z$.

It turns out that every $v_\theta$ has a unique maximizing measure, and that this measure is typically periodic: for Lebesgue almost every $\theta\in\R /\Z$, the $v_\theta$-maximizing measure is supported on a single periodic orbit.
Periodic maximizing measures are also typical in the topological sense: the set
$\{\theta\in \R/\Z:\m_{\max}(v_\theta)\text{ is a periodic singleton}\}$ contains an open dense subset of $\R/\Z$,
and consequently 
$\{v\in V: \m_{\max}(v)\text{ is a periodic singleton}\}$ contains an open dense subset of $V\equiv \R^2$.
In summary, this model problem exhibits \emph{typically periodic optimization}, a phenomenon which has 
subsequently been established for various natural (infinite-dimensional) function spaces $V$ (see \S \ref{typicalsection} for further details).

More can be said about the maximizing measures arising in this specific model problem.
The only \emph{periodic} measures which are maximizing for some degree-one trigonometric polynomial 
are those on which the action of $T$ is combinatorially equivalent to a rational rotation,
while the non-periodic measures which are maximizing for some $v\in V$ correspond to irrational rotations
(their support is a $T$-invariant Cantor set reminiscent of those arising for so-called Denjoy counterexamples in the theory of degree-one circle maps, cf.~e.g.~\cite{veerman}).
More precisely, the maximizing measures for (non-zero) functions in $V$ are \emph{Sturmian} measures:
the Sturmian measure of rotation number $\varrho\in\R/Z$
is 
the
push forward of Lebesgue measure on $X$
under the map $x\mapsto\sum_{n\ge0} \chi_{[1-\varrho,1)}(\{x+n\varrho\})/2^{n+1}$,
where $\{\,\cdot\,\}$ denotes reduction modulo 1.
For example all $T$-invariant measures supported on a periodic orbit of period $<4$ are Sturmian,
though the measure supported on $\{1/5,2/5,3/5,4/5\}\equiv 0011$ is not, and periodic orbits supporting Sturmian measures become increasingly rare as the period grows (see e.g.~\cite{adjr, bousch1, bullettsentenac, jenkinson2, announce, major, morsehedlund} for further details on Sturmian measures and orbits).
Bousch \cite{bousch1} showed that every Sturmian measure arises as the maximizing measure for some $v_\theta$,
and that if $\varrho$ is irrational then $\theta=\theta(\varrho)$ is unique.

The fact that Sturmian measures are precisely the maximizing measures
for this model problem does rely, to an extent, on the particular choice of $f$ and $g$,
though the presence of Sturmian measures is not altogether surprising:
it has subsequently been shown that Sturmian measures arise naturally as maximizing measures
in a variety of similar settings, as will be described in \S \ref{ergodicdominancesection}.

\section{Ergodic optimization as zero temperature thermodynamic formalism}\label{zerotemperaturesection}


Given $(X,T,f)\in\ccc$,
the \emph{pressure} $P(f)=P(T,f)$ is defined as
\begin{equation}\label{pressuredef}
P(f) =\sup_{m\in\m_T} \left(  \int f\, dm + h(m)  \right) \,,
\end{equation}
where $h(m)$ denotes the entropy of $m$.
Any $m\in\m_T$ attaining the supremum
in (\ref{pressuredef}) is called an \emph{equilibrium measure} 
(denoted by $m_f$ if it is unique) for the function $f$ (which
in this context is referred to as a \emph{potential}).
If $f$ is replaced by $tf$ for $t\in\R$, then the entropy term in 
the supremum (\ref{pressuredef}) loses relative importance as $t\to\infty$
(the thermodynamic interpretation of the parameter $t$ 
is as 
an \emph{inverse temperature},
so that letting $t\to\infty$ is referred to as a \emph{zero temperature limit}).
For large values of $t$, an equilibrium measure for $tf$ is almost maximizing for $f$, in that
its integral is close to the maximum ergodic average $\beta(f)$.
More precisely, a number of early authors \cite{coelho, contreraslopesthieullen, conzeguivarch, jenkinson1, jenkinsonetds, parry}
observed,
in various broadly similar settings (with $T$ hyperbolic and $f$ H\"older continuous, so that $m_{tf}$ exists and is unique) 
that the family $(m_{tf})$ has at least one accumulation point $m$
as $t\to\infty$, that
$m$ is an $f$-maximizing measure, and that
$\lim_{t\to\infty} h(m_t) =  h(m)=\max \{ h(\mu) : \mu\in\m_{\max}(f)\}$ (i.e.~any zero temperature accumulation point 
is of maximal entropy among
the set of $f$-maximizing measures). 
Indeed these conclusions are true in wider generality: 
if $X$ is compact, and the entropy map $\mu\mapsto h(\mu)$ is upper semi-continuous\footnote{Upper semi-continuity of entropy
holds if $T$
is expansive (see \cite{walters}), or more 
generally if $T$
admits a finite generating partition 
(see \cite[Cor.~4.2.5]{keller});
in particular this includes all symbolic systems.
Upper semi-continuity is also guaranteed (see \cite{newhouse})
whenever $T$ is a $C^\infty$ map of a compact manifold.},
then every continuous function has at least one equilibrium measure (see \cite[Thm.~9.13(iv)]{walters}), and
it is not hard to establish the following result.

\begin{theorem}\label{basiczerotemperaturetheorem}
{\bf (Zero temperature limits as maximal entropy maximizing measures)}

\noindent
Let $(X,T,f)\in\ccc$ be
such that the entropy map on $\m_T$ is
upper semi-continuous.
For $t\in\R$, if $m_t$ is an equilibrium measure for $tf$
then $(m_t)$ has at least one accumulation point $m\in\m_T$ as $t\to\infty$, and:
\item (i)
$m$ is an $f$-maximizing measure,
\item (ii)
$h(m)=\max \{ h(\mu) : \mu\in\m_{\max}(f)\}$,
\item (iii)
$\lim_{t\to\infty} h(m_t) = h(m)$.
\end{theorem}


In particular,
under the hypotheses of Theorem \ref{basiczerotemperaturetheorem}, if $\m_{\max}(f)=\{m\}$ 
then $m_{t} \to m$ as $t\to\infty$, so
Theorem \ref{genericuniqueness} implies that for typical
$f$ the weak$^*$ limit $\lim_{t\to\infty} m_{t}$ exists,
and is characterised as being the unique $f$-maximizing measure.
A wider investigation of the nature of the set of accumulation points of $(m_t)$, and of whether  $\lim_{t\to\infty} m_t$
\emph{always} exists, was initially focused on the case of $(X,T)$ a subshift of finite type and $f$ locally constant (hypotheses guaranteeing that the unique equilibrium measure $m_{tf}$ is Markov); it was found
\cite{coelho, jenkinsonetds, pollicottsharp1} that limits $\lim_{t\to\infty} m_{tf}$ are not necessarily ergodic,
nor necessarily the evenly weighted
centroid of ergodic maximizing measures of maximal entropy.
In this setting,
the convergence question was resolved by Br\'emont
\cite{bremont1}, who showed\footnote{The paper \cite{bremont1} uses ideas from analytic geometry (semi-algebraic and sub-analytic maps) which are outside the standard toolkit of most ergodic theorists, and despite its 
elegant brevity, the approach of \cite{bremont1} has not subsequently been pursued.}
 that the zero temperature limit \emph{does} always exist, even when $\m_{\max}(f)$
is not a singleton:

\begin{theorem} \label{bremonttheorem} \cite{bremont1} {\bf (Zero temperature convergence for locally constant functions)}

\noindent
For $(X,T)$ a subshift of finite type, and $f:X\to\R$ locally constant,
 $\lim_{t\to\infty} m_{tf}$ exists;
indeed 
$
L_p(X)=\left\{\lim_{t\to\infty} m_{tf}: f\in C(X) \text{ depends on $p$ coordinates}\right\}
$
is finite
for each $p\in\N$.
\end{theorem}

For example, given $(X,T)$ the full shift on two symbols, 
the set  $L_2(X)$ 
 has cardinality 7, and its elements can be listed explicitly (see \cite{bremont1}).
For larger $p$, and for other subshifts of finite type $(X,T)$, the set of possible limits $L_p(X)$
becomes harder to describe.
Progress on this problem
 was made initially by Leplaideur \cite{leplaideur}, then by Chazottes, Gambaudo \& Ugalde \cite{chazottesgambaudougalde} and Garibaldi \& Thieullen \cite{garibaldithieullen2}, using a variety of techniques,
 and can be summarised as follows:

\begin{theorem} {\bf (Description of zero temperature limit for locally constant functions)}

\noindent
If $(X,T)$ is a subshift of finite type, and $f:X\to\R$ is locally constant,
 then $m=\lim_{t\to\infty} m_{tf}$ is concentrated on a certain subshift of finite type $X_f$
 which is itself a finite union of transitive subshifts of finite type.
The finitely many ergodic components $m_i$ of $m=\sum_{i=1}^q \alpha_i m_i$ are unique equilibrium measures of auxiliary potential functions;
these potentials,
and the weights $\alpha_i$,
 can be constructed algorithmically.
\end{theorem}


In the more general setting of Lipschitz functions on subshifts of finite type, the question of whether
zero temperature limits always exist remained open for several years, 
being
finally\footnote{It was noted in \cite{chazotteshochman} that van Enter \& Ruszel \cite{vanenterruszel} had
already given an example of non-convergence in the zero temperature limit, albeit in a somewhat different context:
a nearest neighbour potential model
with the shift map acting on a subset of $(\R/\Z)^\Z$, the significant difference being that the state space
$\R/ \Z$ is non-discrete.}
 settled negatively by Chazottes \& Hochman
\cite{chazotteshochman}:

\begin{theorem}\label{chazotteshochmantheorem}
\cite{chazotteshochman} {\bf (Zero temperature non-convergence)}

\noindent
For $(X,T)$ the full shift on two symbols, there exist Lipschitz functions
$f:X\to\R$ for which 
$\lim_{t\to\infty} m_{tf}$ does not exist.
Indeed such $f$ may be defined as $f(x)=-\text{dist}(x,Y)$, where
$Y\subset X$ is a (carefully constructed) subshift.
\end{theorem}

The flexibility of the approach in \cite{chazotteshochman} allows the full shift in Theorem \ref{chazotteshochmantheorem} to be replaced
by any (one-sided or two-sided) mixing subshift of finite type, 
and allows the construction of subshifts $Y$ such that
the set of accumulation
points of $(m_{tf})$ is e.g.~non-convex, or only containing positive entropy measures,
or not containing ergodic measures.
Bissacot, Garibaldi \& Thieullen \cite{bissacotgaribaldithieullen} have shown that non-convergence in the zero temperature limit can arise for certain
functions on the full 2-shift which take only countably many values, and where the only ergodic maximizing measures are the Dirac measures at the two fixed points.
Yet another approach to non-convergence in the zero temperature limit has been introduced by Coronel \& Rivera-Letelier
\cite{coronelriveraletelier}, partially based on the methods of \cite{vanenterruszel}, establishing a certain persistence of the non-convergence phenomenon:

\begin{theorem} \cite{coronelriveraletelier} {\bf (Persistence of zero temperature non-convergence)}

\noindent
For $(X,T)$ a full shift on a finite alphabet, there exists a Lipschitz function $f_0:X\to\R$, and complementary open subsets $U^+$ and $U^-$ of $X$, such that for any sequence of positive reals $t_i\to\infty$, there is an arbitrarily
small Lipschitz perturbation $f$ of $f_0$ such that the sequence $m_{t_if}$ has an accumulation point whose support
lies in $U^+$, and an accumulation point whose support lies in $U^-$.
\end{theorem}

Temporarily widening our notion of dynamical system to include \emph{higher dimensional} shifts\footnote{Zero temperature non-convergence results for higher dimensional shifts are also proved in \cite{coronelriveraletelier}.}
(i.e.~$G$-actions on $X=F^G$,
where $G=\Z^d$ or $\N^d$ for some integer $d\ge2$, and $F$ is finite),
the following result\footnote{The proof of Theorem \ref{chazotteshochmanlocallyconstanttheorem} 
in \cite{chazotteshochman} relied on work of Hochman
\cite{hochman} establishing that certain one-dimensional subshifts can be simulated in \emph{finite type} 
subshifts of dimension $d=3$; this fact has now been generalised to dimension $d=2$ (see \cite{aubrunsablik, durandromashchenkoshen}), 
suggesting that Theorem \ref{chazotteshochmanlocallyconstanttheorem} is probably valid for all $d\ge 2$
(though certainly not for $d=1$, in view of Theorem \ref{bremonttheorem}).}
  of \cite{chazotteshochman} represents an interesting counterpoint to Theorems
\ref{bremonttheorem} and \ref{chazotteshochmantheorem}:

\begin{theorem}\label{chazotteshochmanlocallyconstanttheorem}
\cite{chazotteshochman} {\bf (Zero temperature non-convergence for locally constant functions on higher dimensional shifts)}

\noindent
For $d\ge 3$, there exist locally constant functions $f$ on $\{0,1\}^{\Z^d}$
such that for every family $(m_{t})_{t>0}$, where $m_{t}$ is an equilibrium measure for $tf$, the limit
$\lim_{t\to\infty} m_{t}$ does not exist.
\end{theorem}


For the case of $(X,T)$ a countable alphabet subshift of finite type, where $X$ is non-compact and
the entropy map $\mu\mapsto h(\mu)$ is not upper semi-continuous,
additional summability and boundedness hypotheses on the locally H\"older function $f:X\to\R$,
together with primitivity assumptions on $X$,
ensure existence and uniqueness of the equilibrium measures $m_{tf}$,
that the family $(m_{tf})$ does in fact have an accumulation point $m$,
and that $h(m) = \lim_{t\to\infty} h(m_{tf}) = \max\{h(\mu): \mu\in \m_{\max}(f)\}$ (see \cite{freirevargas, jenkinsonmauldinurbanski2,morrisjstatphys}), representing an analogue of Theorem \ref{basiczerotemperaturetheorem}.
If in addition $f$ is locally constant,
Kempton \cite{kempton} (see also \cite{freirevargas}) has established the analogue of Theorem \ref{bremonttheorem},
guaranteeing the weak$^*$ convergence of $(m_{tf})$ as $t\to\infty$.
Iommi \& Yayama 
\cite{iommiyayama}
consider almost additive sequences $\mathcal{F}$ of continuous functions defined
on appropriate countable alphabet subshifts of finite type, proving that 
the family of equilibrium measures $(m_{t\mathcal{F}})$ is tight (based on \cite{jenkinsonmauldinurbanski2}),
hence has a weak$^*$ accumulation point, and that any such accumulation point is a maximizing measure for $\mathcal{F}$ (see also \cite{chenzhao, garibaldigomes, sturmanstark, zhao} for general
ergodic optimization in the context of sequences of functions $\mathcal{F}$).

Zero temperature limits have been analysed for certain specific families of functions:
in \cite{jenkinsonetds} for $T$ the doubling map and $f$ 
a degree-one trigonometric polynomial,
in \cite{baravieralopesmengue} 
a specific class of functions defined on the full shift on two symbols and taking countably many values, in
\cite{baravieraleplaideurlopes1} a one-parameter family of functions  defined
on the full shift on three symbols, 
each 
sharing the same two ergodic maximizing measures, and in 
\cite{bclms, lopesmengue2} for the XY model of statistical mechanics.
Connections with large deviation theory have been studied in
\cite{baravieralopesthieullen, lopesmengue1,lopesmohrsouzathieullen},
and the role of the flatness of the potential function has been investigated in \cite{leplaideurjstatphys}.

One source of interest in zero temperature limits of equilibrium measures is \emph{multifractal analysis}, i.e.~the study
of level sets of the form $K_\gamma = \{x\in X: \lim_{n\to\infty} \frac{1}{n}S_nf(x)=\gamma\}$. Each $K_\gamma$ is $T$-invariant, and the \emph{entropy spectrum of Birkhoff averages}, i.e.~the function
$H:[\alpha(f),\beta(f)]\to\R_{\ge0}$ defined by\footnote{The topological entropy $h_{top}(K_\gamma)$ of the (in general non-compact) invariant set $K_\gamma$ is as defined by Bowen \cite{bowen}, or equivalently by
Pesin \& Pitskel' \cite{pesinpitskel}.}
 $H(\gamma)=h_{top}(K_\gamma)$, is in certain (hyperbolic) settings described by the family of equilibrium measures $(m_{tf})_{t\in\R}$, in the sense that $\Gamma:t\mapsto \int f\, dm_{tf}$
is a homeomorphism $\R\to(\alpha(f),\beta(f))$, and (see \cite{lanford}, and e.g.~\cite{bohrrand, hjkps, hp, pesinbook})
\begin{equation*}\label{entropyfn}
H(\gamma)=h(m_{\Gamma^{-1}(\gamma) f}) =\max\left\{h(\mu):\mu\in\m_T, \int f\, d\mu = \gamma\right\}
\quad\text{for all }\gamma\in(\alpha(f),\beta(f))\,.
\end{equation*}
The function $H$ is concave, and extends continuously to the boundary of $[\alpha(f),\beta(f)]$, though the
absence of equilibrium measures $m_{tf}$ with $\int f\, dm_{tf}$ on the boundary prompted investigation of extremal measures (see \cite{coelho,jenkinson1, jenkinsonetds, pollicottsharp1}), and of the (typical) values $H(\alpha(f))$ and $H(\beta(f))$ (see \cite{schmeling}).

Finally, we note that zero temperature limits of equilibrium measures have been studied in a variety of other dynamical settings, including 
Frenkel-Kontorova models
 \cite{anantharaman},
quadratic-like holomorphic maps
 \cite{coronelriveraletelierquadratic},
 multimodal interval maps \cite{iommitodd},
and H\'enon-like maps \cite{takahasi}.

\section{Revelations}\label{revelationsection}

The fundamental problem of ergodic optimization is to say something about maximizing measures.
A most satisfactory resolution is to explicitly identify the $f$-maximizing measure(s) for a given $(X,T,f)\in\ccc$,
though in some cases we may be content with an approximation to an $f$-maximizing measure, or a result asserting that $\m_{\max}(f)$ lies in some particular subset of $\m_T$.
More generally, for a given $(X,T)\in\d$ and a subset $U\subset C(X)$, we may hope to identify a subset
$\n\subset \m_T$ such that $\m_{\max}(f)\subset\n$ for all $f\in U$, or instead
$\m_{\max}(f)\subset\n$ for all $f$ belonging to a large subset of $U$.

In a variety of such settings, it has been noted that a key technical tool is a function we shall 
refer\footnote{The terms \emph{revelation} and \emph{revealed function} are introduced here,
since despite the ubiquity of these concepts there is as yet no established consensus on terminology. 
The function $\psi = \varphi-\varphi\circ T$ we call a revelation has
previously been
referred to as a \emph{sub-coboundary} 
or the solution of
a \emph{sub-cohomology equation},
and the function $\varphi$ in this context has been called a \emph{sub-action},
a \emph{Barabanov function}, a \emph{transfer function}, or a \emph{maximizing function}.
The notion of a revealed function has sometimes gone by the name of a \emph{normal form},  or
in the context of joint spectral radius problems corresponds to a \emph{(maximizing) Barabanov norm}.}
 to as a \emph{revelation}, and an associated result we shall refer to as a \emph{revelation theorem}
 (see \S \ref{revelationtheoremssection}).
First we require the following concept, describing a situation where
the ergodic optimization problem
is easily solved: 

\begin{defn}
Given $(X,T)\in\d$, we say $f\in C(X)$ is \emph{revealed} if its set of maxima
$f^{-1}(\max f)$ contains a compact $T$-invariant set.
\end{defn}


In the (rare) cases when the function $f$ is revealed, it is clear that the maximum ergodic average
$\beta(f)$ equals $\max f$, and that the set of $f$-maximizing measures  is precisely the 
(non-empty) set of $T$-invariant measures whose support is contained in $f^{-1}(\max f)$.

More generally, if we can find $\psi\in C(X)$ satisfying 
\begin{equation}\label{weakcoboundaryequation}
\int \psi\, d\mu=0\quad \text{for all }\mu\in\m_T\,,
\end{equation}
and such that $f+\psi$ is revealed,
then
$\beta(f)=\beta(f+\psi)$ equals $\max (f+\psi)$, 
and $\m_{\max}(f)=\m_{\max}(f+\psi)$
  is precisely the set of $T$-invariant measures whose support is contained in the set $(f+\psi)^{-1}(\max (f+\psi))$.

A natural choice of function $\psi$ satisfying (\ref{weakcoboundaryequation}) is 
a \emph{continuous coboundary}, i.e.~$\psi=\phi-\phi\circ T$ for some $\phi\in C(X)$,
and the ergodic optimization literature has focused mainly 
(though not exclusively, see e.g.~\cite{morris1})
on such $\psi$, since for practical purposes it usually suffices.
This motivates the following definition:

\begin{defn}\label{revelationdefinition}
For $(X,T,f)\in\ccc$, a continuous coboundary $\psi$ is called a \emph{revelation}
if $f+\psi$ is a revealed function, i.e.~if
\begin{equation}\label{maximizingrevelationequation}
(f+\psi)^{-1}(\max (f+\psi))\ \text{contains a compact $T$-invariant set.}
\end{equation}
\end{defn}

Formalising the above discussion, we record the following:

\begin{prop}\label{recordedprop}
If $\psi$ is a revelation for $(X,T,f)\in\ccc$, then
$\beta(f) = \max(f+\psi)$, and
$$\m_{\max}(f) = \m_{\max}(f+\psi) = 
\{\mu\in\m: \text{supp}(\mu)\subset (f+\psi)^{-1}(\max (f+\psi)) \} \neq \emptyset \,.$$
\end{prop}

A consequence of Proposition \ref{recordedprop} is that if $(X,T,f)\in\ccc$ has a revelation then it enjoys the following property (referred to in \cite{bousch2, morris1} as the \emph{subordination principle}): if $\mu\in\m_T$ is $f$-maximizing,
and if the support of $\nu\in\m_T$ is contained in the support of $\mu$, then $\nu$ is also $f$-maximizing.

\begin{example}\label{lengthtwofamily}
\item[\, (a)]
If $T(x)=2x \pmod 1$, 
the function
$f(x) = (2\cos 2\pi x -1)(\sin 2\pi x + 1)$ is not revealed,
but can be written as
$f=g-\psi$
where
$g(x)=2\cos 2\pi x - 1$ is revealed,
and $\psi(x)=\sin 2\pi x - \sin 4\pi x$ is a continuous coboundary, hence a revelation for $f$. 
The unique $f$-maximizing measure is 
therefore the $g$-maximizing measure, namely the Dirac measure
$\delta_0$.
\item[\, (b)]
If $(X,T)$ is the full shift on the alphabet $\{0,1\}$, and the functions $(f_\theta)_{\theta\in\R}$ 
are defined by $f_\theta(x)= f_\theta((x_i)_{i=1}^\infty)=\theta x_1 +x_2 -(\theta+2)x_1 x_2$,
then $f_\theta$ is revealed if and only if $\theta=1$.
For all $c\in\R$, the function $\psi_c(x)=c(x_1-x_2)$ is a coboundary, and
$(f_\theta+\psi_c)(x)=(\theta+c)x_1+(1-c)x_2-(\theta+2)x_1x_2$.
If $\theta> -1$ then $\psi_{(1-\theta)/2}$ is a revelation, and reveals the invariant measure supported on the period-2 orbit to be the unique $f_\theta$-maximizing measure.
If $\theta <-1$ then 
 $\psi_1$ is a revelation, with unique $f_\theta$-maximizing measure
the Dirac measure concentrated on the fixed point $\overline{0}$.
If $\theta=-1$ then $\psi_{(1-\theta)/2}=\psi_1$ is a revelation, and reveals that the $f_{-1}$-maximizing measures
are those whose support is contained
in the golden mean subshift of finite type. 
\end{example}

In this article we have chosen to interpret optimization as maximization, while noting that
the minimizing measures for $f$
are the maximizing measures for $-f$, and
that the minimum ergodic average $\alpha(f) = \min_{\mu\in\m_T} \int f\, d\mu$
is equal to $-\beta(-f)=-\max_{\mu\in\m_T} \int (-f)\, d\mu$.
Occasionally there is interest in simultaneously considering the maximization and minimization problems;
indeed the above discussion 
suggests the possibility of
simultaneously revealing both the minimizing and maximizing measures, by a judicious choice of revelation.
This possibility was considered by Bousch \cite{bousch3}, who showed
(see Theorem \ref{bilateralrevelationtheorem} below)
 that if the $f$-maximizing measures can be revealed, and if the $f$-minimizing measures can be revealed, then indeed it \emph{is} possible to reveal both maximizing and minimizing measures simultaneously.

To make this precise, let us introduce the following terminology. For a given dynamical system $T:X\to X$,
a \emph{revelation} for $f$, in the sense of Definition \ref{revelationdefinition}, will also be called a \emph{maximizing revelation}, while a revelation for $-f$ will be called a \emph{minimizing revelation} for $f$
(i.e.~a minimizing revelation $\psi$ is a continuous coboundary such that
$(f+\psi)^{-1}(\min (f+\psi))$
contains a compact $T$-invariant set).
We say that $\psi$ is a \emph{bilateral revelation} for $f$ if it is both a minimizing revelation and a maximizing revelation. 

\begin{theorem}\label{bilateralrevelationtheorem} \cite{bousch3} {\bf (Bilateralising the maximizing and minimizing revelations)}

\noindent
For $(X,T,f)\in\ccc$, if there exists both a 
minimizing and a maximizing revelation, then there exists a bilateral revelation (i.e.~a continuous coboundary
$\varphi-\varphi\circ T$ with $(f+\varphi-\varphi\circ T)(X)=[\alpha(f),\beta(f)]$).
\end{theorem}

For example, revisiting the family $(f_\theta)_{\theta\in\R}$ from Example \ref{lengthtwofamily}(b), 
if $-3\le \theta\le -2$ 
then 
$\psi_c$ is seen to be a bilateral revelation for $f_\theta$, for all $c\in [-\theta-1,2]$.

\section{Revelation theorems}\label{revelationtheoremssection}

By a \emph{revelation theorem}\footnote{Our terminology is consistent with that of \S \ref{revelationsection}, as again
there is no established consensus on how to describe such theorems: the revelation theorem has been variously called
the \emph{normal form theorem}, the \emph{positive Livsic theorem}, \emph{Ma\~n\'e's lemma},
the \emph{Bousch-Ma\~n\'e cohomology lemma}, and the \emph{Ma\~n\'e-Conze-Guivarc'h lemma}.}
we mean a result of the following kind:

\begin{theorem}\label{modelrevelationtheorem} {\bf (Revelation theorem: model version)}

\noindent 
For a given (type of) dynamical system $(X,T)\in\d$, and a given (type of) function $f\in C(X)$,
there exists a revelation $\varphi-\varphi\circ T$ (i.e.~$f+\varphi-\varphi\circ T$ is a revealed function).
\end{theorem}

In a typical revelation theorem, $(X,T)$ is assumed to enjoy some hyperbolicity, and there is some restriction on the modulus of continuity of $f$.
This is reminiscent of Livsic's Theorem (see \cite{livsic}, or
e.g.~\cite{katokhasselblatt, parrypollicott}),
which 
asserts that if $(X,T)$ is suitably hyperbolic, and $f$ is suitably regular (e.g. H\"older continuous) such that
$\int f\, d\mu=0$ for all $\mu\in\m_T$, then $f$ is a continuous coboundary.
Indeed a Livsic-type theorem can be viewed as a special case of a revelation theorem,
as it follows by
applying an appropriate revelation theorem to both $f$ and $-f$, then invoking Theorem \ref{bilateralrevelationtheorem}.

Revelation theorems date back to the 1990s: Conze \& Guivarc'h proved a version as part of 
\cite{conzeguivarch}, there are parallels with work of Ma\~n\'e on Lagrangian flows \cite{mane1,mane2}, while the first published revelation theorem resembling Theorem \ref{modelrevelationtheorem} was due to 
Savchenko\footnote{Savchenko's 3-page paper contains no discussion of why the revelation theorem is interesting or useful, though does
include some comments on its genesis: the problem had been proposed in Anosov \& Stepin's Moscow dynamical systems seminar in November 1995, and had also
been conjectured by Bill Parry.
Savchenko's proof relies on thermodynamic formalism, 
exhibiting $\varphi$ as a sub-sequential limit of 
$\frac{1}{t} \log h_{t} $, where $h_t$ is the eigenfunction for the dominant eigenvalue
of the Ruelle operator with potential function $tf$.}
 \cite{savchenko}, for $(X,T)$ a subshift of finite type and $f$ H\"older continuous.
Other pioneering papers containing revelation theorems were those of Bousch \cite{bousch1}
and Contreras, Lopes \& Thieullen \cite{contreraslopesthieullen}.

Common features of these early revelation theorems 
are that $f$ is H\"older or Lipschitz, and that $T$ is 
expanding.
The following revelation theorem for expanding maps is 
a particular case of a result in \cite{bousch5} (which is valid for more general \emph{amphidynamical systems}),
and recovers those in \cite{bousch1, contreraslopesthieullen, conzeguivarch, savchenko}:

\begin{theorem} \label{expandinglipschitzrevelationtheorem}
{\bf (Revelation theorem: expanding $T$, Lipschitz $f$)}

\noindent
For expanding $(X,T)\in\d$, 
every Lipschitz function $f:X\to\R$ has a Lipschitz revelation.
\end{theorem}

Since an $\alpha$-H\"older function for the metric $d$ is a Lipschitz function for the 
metric
 $d_\alpha$ defined
by $d_\alpha(x,y)=d(x,y)^\alpha$, we deduce:

\begin{cor} {\bf (Revelation theorem: expanding $T$, H\"older $f$)}

\noindent
For expanding $(X,T)\in\d$,
every 
$\alpha$-H\"older function $f:X\to\R$ 
has an
$\alpha$-H\"older revelation, for all $\alpha\in(0,1]$.
\end{cor}

To prove Theorem \ref{expandinglipschitzrevelationtheorem}, we claim that the function $\varphi$ defined
 by 
 \begin{equation}\label{phiforproof}
 \varphi(x)= \sup_{n\ge1} \sup_{y\in T^{-n}(x)} \left( S_nf(y) - n\beta(f)\right)
 \end{equation}
is such that $\varphi - \varphi\circ T$ 
is a Lipschitz revelation.
Without loss of generality we may assume that $\beta(f)=0$, so that 
(\ref{phiforproof}) becomes 
$$\phi(x)=  \sup_{n\ge1} \sup_{y\in T^{-n}(x)}  S_nf(y)\,.$$
To show that $\phi-\phi\circ T$ is a revelation, we first claim\footnote{Note that the proof of this claim does not require
that $f$ is Lipschitz or that $T$ is expanding.}
 that
\begin{equation}\label{firstclaim}
f+\phi-\phi\circ T \le 0\,,
\end{equation}
and note this immediately implies that $(f+\phi-\phi\circ T)^{-1}(0)$ contains a compact $T$-invariant set, since otherwise there could not be any $m\in\m_T$ satisfying $\int f\, dm=0=\beta(f)$.
To prove (\ref{firstclaim}), note that 
\begin{equation}\label{phiT}
\phi(Tx) = \sup_{n\ge1}\sup_{y\in T^{-n}T(x)} S_nf(y)
\ge \sup_{n\ge1}\sup_{y\in T^{-(n-1)}(x)} S_nf(y)\,,
\end{equation}
because $T^{-(n-1)}(x)\subset T^{-n}(T(x))$.
Now if $y\in T^{-(n-1)}(x)$ then $S_nf(y) = f(x) + S_{n-1}f(y) $ for all $n\ge 1$ (with the usual convention that $S_0f\equiv 0$), so (\ref{phiT}) gives
\begin{equation}\label{phiTagain}
\phi(Tx) 
\ge f(x) + \sup_{n\ge1}\sup_{y\in T^{-(n-1)}(x)} S_{n-1}f(y)\,.
\end{equation}
However, 
\begin{equation}\label{manipulation}
\sup_{n\ge1}\sup_{y\in T^{-(n-1)}(x)} S_{n-1}f(y) 
= \sup_{N\ge0}\sup_{y\in T^{-N}(x)} S_Nf(y)
\ge \sup_{N\ge1}\sup_{y\in T^{-N}(x)} S_Nf(y)
=\phi(x)\,,
\end{equation}
so combining (\ref{phiTagain}) and (\ref{manipulation}) gives
$\phi(Tx) \ge f(x)+\phi(x)$, which is the desired inequality (\ref{firstclaim}).

To complete the proof of Theorem \ref{expandinglipschitzrevelationtheorem}
it remains to show that $\phi$ is Lipschitz, i.e.~that there exists $K>0$ such that for all $x,x'\in X$,
\begin{equation}\label{plusprime}
\phi(x) - \phi(x') \le K d(x,x')\,.
\end{equation}
Given $x,x'\in X$, for any $\varepsilon>0$ there exists $N\ge1$ and $y\in T^{-N}(x)$ such that
\begin{equation}\label{plusplus}
\phi(x)\le S_Nf(y)+\varepsilon\,.
\end{equation}
Since $T$ is expanding we may write $y=T_{i_1}\circ\cdots\circ T_{i_N}(x)$ where the $T_{i_j}$ denote inverse branches of $T$ (i.e.~each $T\circ T_{i_j}$ is the identity map), and we  now define $y'\in X$ by
$y':=T_{i_1}\circ\cdots\circ T_{i_N}(x')$. In particular $y\in T^{-N}(x')$, so
$S_Nf(y') \le \sup_{z\in T^{-N}(x')} S_Nf(z)$, and therefore
\begin{equation}\label{fourplus}
S_Nf(y') \le \sup_{n\ge1} \sup_{z\in T^{-n}(x')} S_nf(z) = \phi(x')\,.
\end{equation}
If $\lambda>1$ is an expanding constant for $T$, i.e.~$d(T(z),T(z'))\ge\lambda d(z,z')$ for all $z,z'$ sufficiently close to each other, then $\gamma=\lambda^{-1}$ is a Lipschitz constant for each of the inverse branches of $T$, so that
if $0\le j\le N-1$
then
$d(T_{i_{j+1}}\circ\cdots\circ T_{i_N}(x),T_{i_{j+1}}\circ\cdots\circ T_{i_N}(x'))\le \gamma^{N-j}d(x,x')$, therefore
$$
f\left( T_{i_{j+1}}\circ\cdots\circ T_{i_N}(x) \right) - f( T_{i_{j+1}}\circ\cdots\circ T_{i_N}(x')) \le Lip(f) \gamma^{N-j}d(x,x') \,,
$$
and hence
\begin{equation}\label{threeplus}
S_Nf(y)-S_Nf(y')
\le \sum_{j=0}^{N-1} Lip(f) \gamma^{N-j}d(x,x')
< \frac{\gamma}{1-\gamma} Lip(f) d(x,x') \,.
\end{equation}
Combining (\ref{plusplus}), (\ref{fourplus}) and
(\ref{threeplus}) then gives
$
\phi(x)-\phi(x') 
<
\frac{\gamma}{1-\gamma} Lip(f) d(x,x') + \varepsilon
$,
but $\varepsilon>0$ was arbitrary, so
$$
\phi(x)-\phi(x') 
\le
\frac{\gamma}{1-\gamma} Lip(f) d(x,x')  \,,
$$
which is the desired Lipschitz condition (\ref{plusprime}) with $K=\frac{\gamma}{1-\gamma}  Lip(f)$, so 
Theorem \ref{expandinglipschitzrevelationtheorem} is proved.

In fact there are various different routes to proving Theorem \ref{expandinglipschitzrevelationtheorem}, stemming
from other possible choices of $\varphi$
 (see e.g.~\cite{garibaldilopes2,garibaldilopesthieullen} for further details), 
notably the choice
$$\varphi(x)=\limsup_{n\to\infty} \sup_{y\in T^{-n}(x)} \left( S_nf(y) - n\beta(f)\right)\,,$$
 which moreover
(see \cite{bousch1, contreraslopesthieullen}) satisfies the functional equation
\begin{equation}\label{eigenequation}
\varphi(x) +\beta(f) = \sup_{y\in T^{-1}(x)} \left( f+\varphi\right)(y)\,.
\end{equation}
Similar functional equations arise in a number of related settings, for example weak KAM theory
\cite{fathi1,fathi2,fathi3,fathi4} and infinite horizon optimal control theory \cite[Thm.~5.2]{chl}.
Indeed (\ref{eigenequation}) can be interpreted as an eigenequation for the operator defined by its righthand side, with 
the maximum ergodic average $\beta(f)$ playing the role of its eigenvalue; the nonlinear operator may be viewed as an analogue of the classical Ruelle transfer operator (see e.g.~\cite{baladibook, parrypollicott, ruellebook}) with respect to the max-plus algebra (in which the $\max$ operation plays the role of addition, and addition plays the role of multiplication, see e.g.~\cite{baravieraleplaideurlopes}).

There is a revelation theorem for maps $T$ which 
satisfy a condition that is weaker than being expanding:
Bousch \cite{bousch2} defined $T:X\to X$ to be \emph{weakly expanding} if its inverse $T^{-1}$ is $1$-Lipschitz
when acting on the set of compact subsets of $X$, equipped with the induced Hausdorff metric
(i.e.~for all $x,y\in X$, there exists $x'\in T^{-1}(y)$ such that $d(x,x')\le d(Tx, Tx')$).
The focus of \cite{bousch2} was on functions which are \emph{Walters} 
(the notion was introduced in \cite{walters2})
for the map $T$: 
for all $\epsilon>0$ there exists $\delta>0$ such that for all $n\in\N$, $x,y\in X$,
if $d(T^ix,T^iy)<\delta$ for $0\le i<n$ then $|S_nf(x)-S_nf(y)|<\epsilon$.

\begin{theorem}\label{waltersrevelation}  \cite{bousch2}
{\bf (Revelation theorem: weakly expanding $T$, Walters $f$)}

\noindent
If $(X,T)\in\d$ is weakly expanding,
then every Walters function $f:X\to\R$ has a revelation.
\end{theorem}



The first revelation theorem in the setting of
\emph{invertible} hyperbolic systems
 appeared in \cite{bousch2}, for maps $T$ satisfying an abstract notion of hyperbolicity dubbed
\emph{weak local product structure}: for all $\epsilon>0$ there exists $\eta>0$ such that if
orbits $(x_i)_{i\le0}$ and $(y_i)_{i\ge0}$ satisfy $d(x_0,y_0)\le \eta$ then there exists an orbit $(z_i)_{i\in\Z}$
with $d(x_i,z_i)\le\epsilon$ for $i\ge0$, and $d(y_i,z_i)\le\epsilon$ for $i\ge0$.

\begin{theorem}\label{bouschwalterstheorem} \cite{bousch2}
{ \bf (Revelation theorem: $T$ with weak local product structure)}

\noindent
If $(X,T)\in\d$ is transitive and has weak local product structure, then every Walters function 
has a revelation.
\end{theorem}

In particular, a transitive Anosov diffeomorphism has weak local product structure,
and in this case
a H\"older continuous function is Walters, so
Theorem \ref{bouschwalterstheorem} implies the existence of a revelation.
The following stronger result
confirms, as suggested by Livsic's Theorem, that in this case the revelation is also H\"older:

\begin{theorem}\label{lopesthiuellenanosovholder} {\bf (Revelation theorem: $T$ Anosov, $f$ H\"older)}

\noindent
If $(X,T)\in\d$ is a transitive Anosov diffeomorphism, and $f$ is $\alpha$-H\"older, then there exists a revelation
$\phi-\phi\circ T$, where $\phi$ is $\alpha$-H\"older.
\end{theorem}

A version of Theorem \ref{lopesthiuellenanosovholder} was proved by Lopes \& Thieullen \cite{lopesthieullen1}, who showed that if $f$ is $\alpha$-H\"older then $\phi$ is $\beta$-H\"older
for some $\beta<\alpha$; the fact that $\phi$ can be chosen with the same H\"older exponent as $f$ was
established by Bousch \cite{bousch5}.

Morris \cite{morris3} considered existence and non-existence of revelations in the context of
circle maps with an indifferent fixed point  (improving on earlier work \cite{branton, souza}). 
Specifically, he considered expanding circle maps of Manneville-Pomeau type $\alpha\in(0,1)$, generalising the Manneville-Pomeau map $x\mapsto x+x^{1+\alpha} \pmod 1$, and proved:

\begin{theorem}\label{morrisintermittenttheorem} \cite{morris3} {\bf (Revelation theorem: $T$ of Manneville-Pomeau type)}

\noindent
If $T$ is an expanding circle map of Manneville-Pomeau type $\alpha\in(0,1)$,
then every H\"older function of exponent $\gamma>\alpha$ has a $(\gamma-\alpha)$-H\"older revelation;
however there exist $\alpha$-H\"older functions which do not have a revelation.
\end{theorem}

The estimate 
on the H\"older exponent of the revelation
in Theorem \ref{morrisintermittenttheorem}
 is sharp: there exist $\gamma$-H\"older functions without any revelation
of H\"older exponent strictly larger than $\gamma-\alpha$ (see \cite{morris3}).
Branco \cite{branco} has considered certain degree-2 circle maps with a super-attracting fixed point, 
proving that if $f$ is  $\alpha$-H\"older, and the super-attracting fixed point is not maximizing,
then there exists an $\alpha$-H\"older revelation. 

If $T:X\to X$ and $f:X\to\R$ are continuous, but $X$ is not compact, there
is no guarantee that $f$-maximizing measures exist: the supremum $\sup_{\mu\in\m_T} \int f\, d\mu$
need not be attained by any $m\in\m_T$.
One way of proving existence of $f$-maximizing measures is to establish a revelation theorem for $f$, 
an approach developed in 
\cite{jenkinsonmauldinurbanskisft, jenkinsonmauldinurbanskinoncompact}, with particular focus on
the case of $(X,T)$ a subshift of finite type on the countable alphabet $\N$. In this setting, if a function
with summable variations is such that its values on a given cylinder set\footnote{Here we use $[i]$ to denote the cylinder set consisting of all sequences $(x_n)_{n=1}^\infty$ such that $x_1=i$.}  are sufficiently larger than its
values `at infinity', then  a revelation exists, and in particular $\m_{\max}(f,T)$ is non-empty.
A prototypical result of this kind
(see \cite{bissacotgaribaldi, jenkinsonmauldinurbanskisft}) 
is:

\begin{theorem}\label{noncompactrevelationtheorem} {\bf (Revelation theorem: non-compact subshift of finite type)}

\noindent
For $(X,T)$ the one-sided full shift on the alphabet $\N$, if $f$ is bounded above, of summable variations, and 
there exists $I\in\N$ with
$\sum_{j=1}^\infty \text{var}_j(f) < \inf f|_{[I]} - \sup f|_{[i]}$,
for all sufficiently large $i$, then $f$ has a revelation, and in particular has a maximizing measure.
\end{theorem}

The inequality in Theorem \ref{noncompactrevelationtheorem} 
clearly holds whenever
$\sup f|_{[i]}\to -\infty$, which in particular is the case if $f$ satisfies the summability condition
$\sum_{i=1}^\infty e^{\sup f|_{[i]}} <\infty$ familiar from thermodynamic formalism (see e.g.~\cite{gdms}).
Note that
\cite{bissacotfreire} provides alternative criteria guaranteeing existence
of a maximizing measure for certain functions $f$ defined on irreducible countable alphabet subshifts of finite type:
the approach is more direct than in \cite{jenkinsonmauldinurbanskisft}, and while it does not prove the existence of a revelation, it does establish the 
subordination principle.


Going beyond the setting of discrete dynamics,
there has been some work on revelation theorems for flows\footnote{Note that although the majority of work on ergodic optimization has been placed in the setting of discrete time,
there have been various developments in the context of flows (see \cite{bousch3, lopesrosasruggiero, lopesthieullen2, lopesthieullen3, pollicottsharp2, yanghuntott}).}:
for $T^t$ a smooth Anosov flow without fixed points, and $f$ H\"older continuous, 
there exists a H\"older 
 function $\phi$ satisfying
$\int_0^s f(T^t(x))\, dt + \phi(x)-\phi(T^s(x)) \le s\beta(f)
$
for all $x\in X$, $s\in\R^+$
(see \cite{lopesthieullen2, pollicottsharp2}),
which moreover is smooth in the flow direction (see \cite{lopesthieullen2}).
An analogous result holds for certain expansive non-Anosov geodesic flows, see
\cite{lopesrosasruggiero}.

\section{Typically periodic optimization (TPO)}\label{typicalsection}

Given a dynamical system $(X,T)\in\d$ of a particular kind (e.g.~enjoying some appropriate hyperbolicity), 
we wish to establish properties of \emph{typical} maximizing measures:
for a
topological vector space $V$ of real-valued functions on $X$,
we aim to show 
there 
exists
$V'\subset V$ which is topologically large 
(e.g.~containing 
an open dense subset of $V$)
such that all $f\in V'$ have maximizing measure(s) with a certain specified property.
The specified property 
we have in mind
is that the maximizing measure be \emph{periodic},
though first we note that a weaker property follows as
a simple consequence of \S \ref{revelationtheoremssection}
(where for definiteness $(X,T)$ is assumed to be expanding or Anosov, and $V=Lip$, so Theorems 
\ref{expandinglipschitzrevelationtheorem} and \ref{lopesthiuellenanosovholder} can be used):

\begin{theorem} \label{tonfs}  {\bf (Typical maximizing measures are not fully supported)}

\noindent
Suppose $(X,T)\in\d$ is either expanding or Anosov, and is transitive but not reduced to a single periodic orbit.
The open dense set $Lip'$, defined as the 
complement in $Lip$ of the closed subspace $\{c+\varphi-\varphi\circ T:c\in\R, \varphi\in Lip\}$,
is such that if $f\in Lip'$ then no
$f$-maximizing measure is fully supported.
\end{theorem}

The possibility of typical maximizing measures being \emph{periodic} was suggested by
the early work on ergodic optimization for finite-dimensional spaces of functions, as described in
\S \ref{modelsection}. 
We state this below as the (purposefully imprecise) Conjecture \ref{metatpo}, but first require some notation.

\begin{defn}
For $(X,T)\in\d$, and $V$ a Banach space consisting of certain continuous real-valued functions on $X$, define
$
V_{Per}
$
to be the set of those $f\in V$
such that $\m_{\max}(f)$
contains at least one measure supported on a single periodic orbit.
\end{defn}

\begin{conjecture}\label{metatpo}
{\bf (Typically Periodic Optimization (TPO) Conjecture)}

\noindent
If $(X,T)\in\d$ is a suitably hyperbolic dynamical system, and $V$ is a Banach space consisting of suitably regular continuous functions, then $V_{Per}$ contains an open dense subset of $V$.
\end{conjecture}

The earliest published paper containing specific articulations of 
the TPO Conjecture 
was that of Yuan \& Hunt
\cite{yuanhunt}, where $(X,T)$ was assumed to be either an expanding map or an Axiom A diffeomorphism;
the analogue of the TPO Conjecture 
was conjectured \cite[Conj.~1.1]{yuanhunt} for $V$ a space of smooth (e.g.~$C^1$) functions on $X$, though the case $V=Lip$ was discussed in more detail.
In subsequent years this case $V=Lip$ became a focus of attention among workers in ergodic optimization,
culminating in its resolution 
(see Theorem \ref{lipschitztpo} below)
by Contreras \cite{contreras}, building on work of \cite{morris2,quassiefken,yuanhunt}.

The first proved (infinite-dimensional) version of the TPO Conjecture 
was due to
Contreras, Lopes \& Thieullen \cite{contreraslopesthieullen}, in a paper prepared at around the same time as
\cite{yuanhunt}.
In the context of (smooth) expanding maps (on the circle),
they noted a significant consequence of the revelation theorem (a version of which they proved 
\cite[Thm.~9]{contreraslopesthieullen}), which would also be exploited by subsequent authors: 
\emph{if} it is known that $(X,T,f)$ has a revelation $\psi$ for every $f\in V$, then the revealed function $f+\psi$ is usually more amenable to analysis, and in particular it may be possible to exhibit small perturbations of $f+\psi$ which lie in the desired set $V_{Per}$,
and thus deduce that $f$ itself can be approximated by members of $V_{Per}$.
Their choice of $V=V^\alpha$ was as a closed subspace of $H^\alpha$, the space of $\alpha$-H\"older functions on $X$; the space $V^\alpha$ is defined to consist of the closure (in $H^\alpha$) of those functions which are actually \emph{better} than $\alpha$-H\"older, i.e.~they are $\beta$-H\"older for some $\beta > \alpha$ (so $V^\alpha$ is defined for $\alpha<1$, but undefined in the Lipschitz class $\alpha=1$).
The superior approximation properties enjoyed by $V^\alpha$ yield:

\begin{theorem}\label{cltlittleholder} \cite{contreraslopesthieullen} {\bf (TPO on a proper closed subspace of H\"older functions)}

\noindent
For $\alpha\in (0,1)$ and $(X,T)$ a circle expanding map, 
 $V^\alpha_{Per}$ contains an open dense subset of $V^\alpha$.
\end{theorem}

 Bousch \cite[p.~305]{bousch2}
was able to use his revelation theorem (Theorem \ref{waltersrevelation})
for the set $W$ of Walters functions in the particular context of the one-sided
full shift, where $W$ can be given the structure of a Banach space,  
to prove the following:

\begin{theorem}\label{bouschwalterstpo}  \cite{bousch2} {\bf (TPO for Walters functions on a full shift)}

\noindent
For $(X,T)$ a full shift, $W_{Per}$ contains an open dense subset of $W$.
\end{theorem}

An important ingredient in the proof of Theorem \ref{bouschwalterstpo} is that locally constant functions are
dense in $W$, and that for such functions $f$
the set $\m_{\max}(f)$ is 
stable under perturbation (indeed $\m_{\max}(f)$
 is the set of all invariant measures supported by some subshift of finite type,
and such subshifts always contain at least one periodic orbit).

Quas \& Siefken \cite{quassiefken} also considered the setting of ($X,T)$ a one-sided full shift, 
and spaces of functions which are Lipschitz with respect to non-standard metrics on $X$:
for a sequence $A=(A_n)_{n=1}^\infty$ with $A_n\searrow0$,
define a metric $d_A$ on $X$ by $d_A(x,y)=A_n$ if $x$ and $y$ first differ in the $n$th position
(i.e.~$x_i=y_i$ for $1\le i<n$, and $x_n\neq y_n$), and let $Lip(A)$
denote the space of functions on $X$ which are Lipschitz with respect to $d_A$, equipped with the
induced Lipschitz norm. 
Quas \& Siefken required the additional condition
$\lim_{n\to0} A_{n+1}/A_n =0$
(in which case members of $Lip(A)$ are referred to as \emph{super-continuous} functions) and proved:

\begin{theorem}\label{quassiefkentheorem} \cite{quassiefken} {\bf (TPO for super-continuous functions)}

\noindent
For $(X,T)$ a full shift, if $\lim_{n\to0} \frac{A_{n+1}}{A_n} =0$
then $Lip(A)_{Per}$ contains an open dense subset of $Lip(A)$.
\end{theorem}

In the same context of super-continuous functions on a one-sided full shift $(X,T)$,
Bochi \& Zhang \cite{bochizhang} 
found a more restrictive condition on the sequence $A$ which suffices to guarantee that $Lip(A)_{Per}$ is a 
\emph{prevalent}\footnote{Prevalence is a probabilistic notion of typicalness, 
introduced by Hunt, Sauer \& Yorke \cite{huntsaueryorke},
and in finite dimensional spaces coincides with the 
property of being of full Lebesgue measure.
Specifically, for $V$ a complete metrizable topological vector space, a Borel set $S\subset V$ is called \emph{shy}
if there exists a compactly supported measure which gives mass zero to every translate of $S$, and a \emph{prevalent} set in $V$ is defined to be one whose complement is shy.} 
subset of $Lip(A)$:

\begin{theorem} \cite{bochizhang} \label{bochizhangtheorem} {(\bf Prevalent periodic optimization})
For $(X,T)$ the one-sided full shift on two symbols, if $\frac{A_{n+1}}{A_n} = O(2^{-2^{n+2}})$ as $n\to\infty$
then  $Lip(A)_{Per}$ is
a prevalent subset of $Lip(A)$.
\end{theorem}

The proof in \cite{bochizhang} uses Haar wavelets to 
reduce the problem
to a finite-dimensional one with a graph-theoretic reformulation as a maximum cycle mean problem.
Since the hard part of proving Theorem \ref{quassiefkentheorem} is to show that $Lip(A)_{Per}$ contains a dense
subset of $Lip(A)$, and any prevalent subset is dense, we note that Theorem \ref{bochizhangtheorem}
constitutes a strengthening of Theorem \ref{quassiefkentheorem} in the case that
$\frac{A_{n+1}}{A_n} = O(2^{-2^{n+2}})$ as $n\to\infty$.

Prior to Contreras' proof of Theorem \ref{lipschitztpo} below, a number of authors 
(notably \cite{bouschsmf, bremontflowers,yuanhunt})
had considered
the case $V=Lip$ in the TPO Conjecture, 
and established partial and complementary results.
The first of these was due to
Yuan \& Hunt \cite{yuanhunt}:

\begin{theorem}  \cite{yuanhunt} \label{yuanhunttheorem}
{\bf (Non-periodic measures are not robustly optimizing)}

\noindent
Let $(X,T)\in\d$ be an expanding map.
If $f\in Lip$ has a non-periodic maximizing measure $\mu$, then there exists $g\in Lip$,
arbitrarily close to $f$ in the Lipschitz topology,
such that $\mu$ is not $g$-maximizing.
\end{theorem}

Bousch \cite{bouschsmf}
gave an alternative proof of Theorem \ref{yuanhunttheorem}, 
 in the more general setting of amphidynamical systems, making explicit the role of revelations,
 and
 quantifying 
 the phenomenon of periodic orbits of low period being more stably 
 maximizing than those of high period: if $f\in Lip$ has a periodic maximizing measure $\mu$ of (large) period $N$,
 then there exist $O(1/N)$-perturbations of $f$ 
 in the Lipschitz norm for which $\mu$ is no longer maximizing. 
 More precisely:
 
 \begin{prop}\label{bouschbound} \cite{bouschsmf} {\bf (A bound on orbit-locking for Lipschitz functions)}
 
 \noindent
 If $(X,T)\in\d$ is expanding, then there exists $K_T>0$ such that if $f\in\text{Lip}$ has $\mu$ as an $f$-maximizing
 measure, where $\mu$ is supported on a periodic orbit of period $N>K_T$, then there exists $g\in\text{Lip}$
 with $\text{Lip}(f-g)\le (\frac{N}{K_T} -1)^{-1}$ such that $\mu$ is not $g$-maximizing.
 \end{prop} 

The constant $K_T$ in Proposition \ref{bouschbound} can be chosen as
 $K_T = 6C_TL_T$, where $L_T$ is a Lipschitz constant for $T$, and $C_T$ is such that
$\text{Lip}(\phi) \le C_T \text{Lip}(f)$ whenever $\phi-\phi\circ T$ is a revelation for a Lipschitz function $f$,
so for example
in the particular case of $T(x)=2x \pmod 1$ on the circle, we may take $L_T=2$ and
 $C_T=1$ (see \cite[Lem.~B]{bousch1}), so
for any $f\in Lip$ whose maximizing measure $\mu$ is supported on an orbit of period $N>12$, there 
exists $g\in Lip$ with $\text{Lip}(f-g) \le 12/(N-12)$ such that $\mu$ is not $g$-maximizing.

A proof of the TPO Conjecture 
in the important case $V=Lip$ was given by Contreras: 

\begin{theorem}\label{lipschitztpo} \cite{contreras} {\bf (TPO for Lipschitz functions)}

\noindent
For $(X,T)\in\d$ an expanding map,
$Lip_{Per}$ contains an open dense subset of $Lip$.
\end{theorem}

To sketch\footnote{This sketch follows the exposition of Bousch \cite{bouschaftercontreras}.} a proof of Theorem \ref{lipschitztpo}, we first note that if $\mu$ is any periodic orbit measure, 
it is relatively easily shown that
$\{f\in Lip: \mu\text{ is } f\text{-maximizing}\}$ is a closed set with non-empty interior, so 
it suffices to show that
$Lip_{Per}$
is dense in $Lip$.
Let us say that $\mu\in\m_T$ is a \emph{Yuan-Hunt measure} if for all $x\in\text{supp}(\mu)$,
$Q>0$, there exist integers $m,p\ge0$ such that\footnote{Condition (\ref{classonedefn}) was introduced by Yuan \& Hunt \cite[p.~1217]{yuanhunt}, who called it the Class I condition. (A Class II condition was also introduced in \cite{yuanhunt}, in terms of approximability by periodic orbit measures, which stimulated related work in \cite{bressaudquas, colliermorris}).}
\begin{equation}\label{classonedefn}
\min\{d(T^ix,T^jx): m\le i,j\le m+p, 0<|i-j|<p\}
> Q d(T^{m+p}x,T^mx)\,.
\end{equation}
If $\m_{YH}$ denotes the set of Yuan-Hunt measures, and $$Lip_{YH}=\{f\in Lip: \m_{\max}(f)\cap \m_{YH}\neq\emptyset\}\,,$$
then clearly every invariant measure supported on a periodic orbit lies in $\m_{YH}$, 
and thus $Lip_{Per}\subset Lip_{YH}$.
Yuan \& Hunt proved \cite[Lem.~4.10]{yuanhunt} that $Lip_{Per}$ is dense in $Lip_{YH}$, so to prove
Theorem \ref{lipschitztpo} it suffices to show that $Lip_{YH}$ is dense in $Lip$.
Contreras \cite{contreras} showed, by estimating the lengths of pseudo-orbits, that if $\mu\in\m_T \setminus \m_{YH}$
then $\mu$ has strictly positive entropy.
However, a result of Morris \cite{morris2} asserts that the set of Lipschitz functions with a positive entropy maximizing measure is of first category; it follows that $Lip\setminus Lip_{YH}$ is of first category, and therefore
$Lip_{YH}$ is dense in $Lip$, as required.

\section{Other typical properties of maximizing measures}\label{othertypicalsection}

For a suitably hyperbolic dynamical system $(X,T)$, the fact that typical properties of maximizing
measures in $C(X)$ are rather different from those of more regular continuous
functions discussed in \S \ref{typicalsection} is illustrated by the following result:

\begin{theorem}\label{conttypical} \cite{morris4} {\bf (Typical maximizing measures for continuous functions)}

\noindent
For $(X,T)\in\d$ either expanding or Anosov, 
and transitive but not reduced to a single fixed point,
there is a residual subset $C'\in C(X)$ such that if $f\in C'$
then $\m_{\max}(f)$ is a singleton containing a measure which is fully supported,
has zero entropy, and is not strongly mixing. 
\end{theorem}

Parts of Theorem \ref{conttypical} had been proved elsewhere
(see \cite{bouschjenkinson} for
the fact that typical maximizing measures are fully supported,
and \cite{bremontcras} for
the fact that typical maximizing measures have zero entropy),
while Morris's proof in \cite{morris4}
was a natural consequence of his following more abstract result (together
with results of Sigmund \cite{sigmund} on residual subsets of the set of invariant measures):

\begin{theorem} \cite{morris4} {\bf (Maximizing measures inherit typical properties from $\m_T$)}

\noindent
Suppose $(X,T)\in\d$ is such that the set of ergodic measures is weak$^*$ dense in $\m_T$.
Then for typical continuous functions $f$, the $f$-maximizing measure inherits any property which is
typical in $\m_T$.
More precisely, if $\m'$ is a residual subset of $\m_T$, then
$\{f\in C(X):\m_{\max}(f)\subset \m'\}$ is a residual subset of $C(X)$. 
\end{theorem}

The most surprising aspect of Theorem \ref{conttypical}
is that in $C(X)$ a typical maximizing measure is \emph{fully supported}.
Not only does this contrast with Theorem \ref{tonfs} and 
the typical \emph{periodic} optimization results
of \S \ref{typicalsection},
but it also contrasts with intuition;
indeed an open problem is to exhibit constructively 
a continuous function $f:X\to\R$, and
an expanding or Anosov dynamical system
$(X,T)\in\d$, such that the unique $f$-maximizing measure is fully supported.  
Clearly such a unique maximizing measure must be ergodic (since the set
$\m_{\max}(f)$ is convex, and ergodic
maximizing measures are precisely its extreme points), though it turns out 
that this is the only restriction: 

\begin{theorem}\label{uniqmaxthm} \cite{jenkinsoneem} {\bf (Every ergodic measure is uniquely maximizing)}

\noindent
If $(X,T)\in\d$
then for any ergodic $\mu\in\m$, there exists $f\in C(X)$ such that $\mu$ is the unique $f$-maximizing measure.
\end{theorem}

The above results, and those
of \S \ref{typicalsection}, 
involve fixing the dynamical system
$(X,T)\in\d$, and enquiring about typical properties of $f$-maximizing measures for $f$ lying
in some Banach space $V$.
More generally, one might view the triple $(X,T,f)\in\ccc$ as varying in some given topological space,
and again enquire about typical properties of $\m_{\max}(X,T,f)$; alternatively we may fix the
(compact metric) state space $X$, and view pairs $(T,f)$ as elements of some topological space ${\mathcal P}$,
and again enquire about typical properties of $\m_{\max}(T,f)$.

The existing literature has focused on three versions of this general problem, the first (and most studied)
of which is the case ${\mathcal P}=\{T\}\times V$ described
previously.
 A second case is when $T$ varies within some topological space ${\mathcal T}$ of maps, and the function
 $f=f_T$ varies non-trivially with $T$. In fact the limited work on this second case has focused
 (see \cite{contreraslopesthieullen,jenkinsonmorris, moritatokunaga}) on
 \emph{Lyapunov maximizing measures} for expanding maps $T$,
 i.e.~where $f_T=\log |T'|$, so that
 $\max_{\mu\in\m_T} \int \log |T'|\, d\mu$ is the \emph{maximum Lyapunov exponent}.
 If ${\mathcal T}$ is the space of $C^1$ expanding maps then a typical Lyapunov maximizing measure is shown
 (see \cite{jenkinsonmorris} for $X$ the circle, and \cite{moritatokunaga} for $X$ a more general manifold)  to
 be unique, fully supported, and of zero entropy (this can be considered an analogue of Theorem
 \ref{conttypical}); by contrast, in the context of $C^{1+\alpha}$ expanding maps an analogue
 of Theorem \ref{cltlittleholder} is established (see \cite{contreraslopesthieullen}), and indeed
 it is likely that other results in \S \ref{typicalsection} have natural analogues in the context of Lyapunov maximizing measures.
 
 The third version of the general problem involves fixing the continuous function $f:X\to\R$
 (with possibly additional hypotheses on $f$) and varying the map $T$ within some topological space ${\mathcal T}$.
For example, if $X$ is a compact connected manifold of dimension greater than 2, 
$T$ is varied in the space $\mathcal{T}$
 of homeomorphisms of $X$,
and $f\in C(X)$ is considered fixed,  
it can be shown (see \cite{talzanatafund}) that there is a \emph{dense} subset $\mathcal{T}' \subset \mathcal{T}$
such that $\m_{\max}(X,T,f)$ contains a periodic orbit measure.
However, provided $f$ is non-constant when restricted to any non-empty open subset, it turns out that periodic maximization is \emph{not} typical:

\begin{theorem}  \label{talzanataresidualtheorem} \cite{talzanatadcds} {\bf (Typical optimization is not periodic)}

\noindent
Let $X$ be a compact  connected Riemannian manifold of dimension at least 2.
If $f\in C(X)$ is non-constant when restricted to any non-empty open subset, then there is a
residual subset  $\mathcal{T}_f' \subset \text{Homeo}(X)$ 
such that for every $T\in \mathcal{T}_f'$,
the set $\m_{\max}(X,T,f)$ contains no periodic orbit measures.
\end{theorem}

Analogous results have been established for 
the space $\mathcal T$ of endomorphisms (i.e.~continuous surjections):
for example in \cite{batistagonschorowskital} (see also \cite{talzanatanonlinearity}) 
it is shown that for any
compact Riemannian manifold $X$, and any $f\in C(X)$, there is a dense subset $\mathcal{T}'\subset \mathcal{T}$
such that $\m_{\max}(X,T,f)$ contains a periodic orbit measure for all $T\in\mathcal{T}'$;
however if $X$ is the circle it is known 
(see \cite{talzanatadcds})
that such a $\mathcal{T}'$ is meagre unless the function $f$
is monotone on some sub-interval (in which case $\mathcal{T}'$ has non-empty interior).

\section{Sturmian optimization and ergodic dominance}\label{ergodicdominancesection}

In ergodic optimization, Sturmian measures were first observed in the context of the model problem described in \S \ref{modelsection}, as maximizing measures for functions of the form
$v_\theta(x)=\cos 2\pi(x-\theta)$, with underlying dynamical system $T(x)=2x \pmod 1$
on the circle $\R/\Z$.
As well as their definition in terms of rotations (see \S \ref{modelsection}), Sturmian measures can be characterised as precisely those $T$-invariant probability measures whose support is contained in a sub-interval of the form $[\gamma,\gamma+1/2]$ (i.e.~a closed semi-circle), see e.g.~\cite{bouschmairesse, bullettsentenac}.
In other words, the family of Sturmian measures can be defined as the maximizing measures  for the family of characteristic functions $\chi_{[\gamma,\gamma+1/2]}$, $\gamma\in\R/\Z$.
In view of this definition, 
it is perhaps not so surprising that Sturmian measures 
arise
as maximizing measures for certain naturally occurring functions,
and indeed they have subsequently been identified as maximizing measures for functions other than the family $v_\theta$.
One such example (see \cite{adjr}) is the family of functions
$u_\theta(x)= - d(x,\theta)$, $\theta\in\R/\Z$, where $d$ is the usual distance function on the circle.
As for the family of functions $v_\theta$, the Sturmian measures are \emph{precisely} the maximizing measures for the functions $u_\theta$, with each non-periodic Sturmian measure 
being maximizing for a single function $u_\theta$, but periodic Sturmian measures 
corresponding to a positive length closed interval of parameters $\theta$.

Moving beyond finite dimensional families of functions, 
there exist infinite dimensional function cones where Sturmian measures are guaranteed to be maximizing;
by a function cone we mean a set $K$ of functions on $X$ which is closed under addition
(i.e.~$K+K\subset K$) and multiplication by non-negative reals (i.e.~$\R_{\ge0} K\subset K$).
If $X=[0,1]$ then the set of concave real-valued functions on $X$ is a cone, and
if $T$ is the doubling map\footnote{Although the doubling map on $[0,1]$ is not continuous, 
its set of invariant probability measures is nevertheless weak$^*$ compact, so $\m_{\max}(f)\neq\emptyset$ for all continuous $f$.}
 on $[0,1]$, with $T(1)=1$ and $T(x)=2x \pmod 1$ for $x<1$, then:

\begin{theorem}\label{sturmianmaximizingconcave} \cite{announce,major} {\bf (Sturmian maximizing measures for concave functions)}

\noindent
For the doubling map $T:[0,1]\to[0,1]$, if $f:[0,1]\to\R$ is concave then it has a Sturmian maximizing measure.
If $f$ is strictly concave then its maximizing measure is unique and Sturmian.
\end{theorem} 

The set of increasing functions on $[0,1]$ is also a function cone. For the doubling map on $[0,1]$, the Dirac measure $\delta_1$ is clearly $f$-maximizing for every increasing function $f:[0,1]\to\R$.
This simple fact has
a more surprising generalisation:
if\footnote{For $\beta>2$ there is a slightly different
version of Theorem \ref{betamaxtheorem} for analogous maps $T_\beta$ (see \cite{aj} for details).}
 $\beta\in(1,2)$, and $T_\beta:[0,1]\to[0,1]$ is given by $T_\beta(x)=\beta x$ on $[0,1/\beta]$
and $T_\beta (x)=\beta x-1$ on $(1/\beta,1]$, then for certain $\beta$ 
(e.g.~the golden mean $\beta=(1+\sqrt{5})/2$)
there exists a single $T_\beta$-invariant
probability measure $\mu_\beta$ which is simultaneously maximizing for all increasing functions on $[0,1]$, and
in this case $\mu_\beta$ is once again 
Sturmian\footnote{Sturmian measures can, as before, be defined in terms of circle rotations; alternatively, 
in this context they are characterised by having support contained in a closed interval of length $1/\beta$.}:

\begin{theorem}\label{betamaxtheorem} \cite{aj} {\bf (Sturmian maximizing measure for all increasing functions)}

\noindent
For the map $T_\beta:[0,1]\to[0,1]$, if $\beta\in(1,2)$ is the dominant root of 
$x^{ap+1}-\sum_{i=0}^p x^{ia}$ for some integers $a,p\ge1$,
then the point 1 generates a Sturmian periodic orbit, and the $T_\beta$-invariant measure on this orbit is
$f$-maximizing for every increasing function $f:[0,1]\to\R$.
\end{theorem}


For any cone $K$ such that
$K-K$ is dense in $C(X)$,
 a partial order
$\prec$ on Borel probability measures arises by declaring that $\mu\prec\nu$ if and only if 
$\int f\, d\mu \le \int f\, d\nu$ for all $f\in K$.
Both the cone of increasing functions and the cone of concave functions enjoy this property, and for these cases
the associated partial order is known as a \emph{stochastic dominance} order
(see e.g.~\cite{bawa, kertzrosler, marshallolkin, shakedshanthikumar}).
We therefore use the term \emph{ergodic dominance} to refer to the study of the partially ordered set $(\m_T, \prec)$,
and the identification of maximal and minimal elements in $(\m_T, K)$
 may be viewed as a generalisation of ergodic optimization.
For $K$ the cone of increasing functions, ergodic dominance in the context of the full shift on two symbols has been investigated in \cite{anagnostopoulou1,anagnostopoulou2}.

For $K$ the cone of concave functions, ergodic dominance has been studied in \cite{jenkinsonsteel} for orientation-reversing expanding maps $T:[0,1]\to[0,1]$,
and in \cite{steel} for certain unimodal maps.
A necessary condition for the comparability of two measures $\mu,\nu$ is that their barycentres coincide, i.e.~$\int x \, d\mu(x)=\int x \, d\nu(x)$, so if $T$ is the doubling map then
$(\m_T,\prec)$ cannot have a maximum element,
though each of the sets $\m_{T,\varrho}=\{\mu\in\m_T:\int x\, d\mu(x)=\varrho\}$ 
does turn out to have such an element:

\begin{theorem} \label{sturmianmaximumelements}  \cite{announce, major} {\bf (Sturmian measures as maximum elements in each $\m_{T,\varrho}$)}

\noindent
If $T:[0,1]\to[0,1]$ is the doubling map, and $\prec$ is the partial order induced by the cone of concave functions on $[0,1]$, then the Sturmian measure 
of rotation number $\varrho$ is the maximum element in $\m_{T,\varrho}$, for all $\varrho\in[0,1]$.
\end{theorem}

In fact Theorem 
\ref{sturmianmaximizingconcave} can be viewed as one of several corollaries to
Theorem \ref{sturmianmaximumelements};
others are that Sturmian measures have strictly smallest variance around their means, and that Sturmian periodic orbits have larger geometric mean than any other periodic orbits with the same arithmetic mean (see \cite{announce, major, balancedmajor} for further details).

Underlying the results in this section is an idea of Bousch  \cite{bousch1}, which provides an approach to proving the Sturmian nature of maximizing measures: the existence of a revelation is guaranteed by
Theorem \ref{expandinglipschitzrevelationtheorem},
and it is potentially feasible to show that the corresponding revealed function takes its maximum value on a closed interval
of length $1/\beta$, in which case the maximizing measure is Sturmian.
This approach can also be used for more general expanding maps $T$ (see e.g.~\cite{jenkinsonhitting, jenkinsonpams, jenkinsonpollicott})
where the closed interval in question has the property that $T$ is injective when restricted to its interior,
and for certain generalisations of Sturmian measures (see e.g.~\cite{bremontflowers, bremont2, hj}).

The article \cite{jenkinsonpollicott} treats a problem concerning the joint spectral radius of pairs of matrices (see e.g.~\cite{dl, jungers, lagariaswang, rotastrang} for background to this area),
which is reformulated as an ergodic optimization problem involving a one-parameter family of expanding maps, and a one-parameter family of functions, whose maximizing measures turn out to be precisely the family of Sturmian measures.
The role of Sturmian measures (or orbits) in this context had previously been noted in
\cite{btv, bouschmairesse, hmst, morrissidorov}.
More generally, 
joint spectral radius problems have a number of parallels with ergodic optimization, and the two fields enjoy a fruitful interaction, see e.g.~\cite{bochimorris, bochirams, bouschmairesse, morrisbarabanov, morris5, morrisbergerwang, morrisplms}.

\end{document}